\documentclass[a4paper,12pt]{article}

\usepackage{amsmath}
\usepackage{amssymb}
\usepackage{url}

\ifx\pdfoutput\undefined
\usepackage{graphicx}
\usepackage{epsfig}
\newcommand{\altgraphicsWidth}[3]{\epsfig{#1,width=#3}}
\else
\usepackage[pdftex]{graphicx}
\newcommand{\altgraphicsWidth}[3]{\includegraphics[width=#3]{#2}}
\fi

\newtheorem{thm}{Theorem}[section]
\newtheorem{lem}[thm]{Lemma}
\newtheorem{prp}[thm]{Proposition}
\newtheorem{cor}[thm]{Corollary}

\newtheorem{con}[thm]{Conjecture}
\newenvironment{proof}{{\it Proof}.\ }{\hfill$\square$\par\medskip}

\def\cC{{\mathcal C}}

\def\cF{{\mathcal F}}

\def\cZ{{\mathcal Z}}

\newcommand\CC{{\mathbb C}}
\newcommand\FF{{\mathbb F}}

\newcommand\RR{{\mathbb R}}
\newcommand\Sph{{\mathbb S}}
\newcommand\ZZ{{\mathbb Z}}

\newcommand\SetOf[2]{\bigl\{#1\,\bigm|\,#2\bigr\}}
\newcommand\id{\operatorname{id}}

\newcommand\norm[1]{|\!|#1|\!|}

\newcommand\GF[1]{\FF_{#1}}

\def\Sym{\operatorname{Sym}}

\def\SetOf#1#2{\left\{\left.#1\vphantom{#2}\ \right|\ #2\vphantom{#1}\right\}}

\def\to{\rightarrow}

\def\st{\operatorname{st}}
\def\lk{\operatorname{lk}}
\def\sd{\operatorname{sd}}

\def\nothing{\rule{0mm}{1ex}}

\title{Projectivities in Simplicial Complexes and Colorings of Simple Polytopes}
\author{Michael Joswig\thanks{Supported by Deutsche Forschungsgemeinschaft, Sonderforschungsbereich~288
    ``Differentialgeometrie und Quantenphysik''.}}
\date{\today}

\def\Sym{\operatorname{Sym}}

\begin{document}
\maketitle

\begin{abstract}
  For each strongly connected finite-dimensional (pure) simplicial complex~$\Delta$ we construct a finite
  group~$\Pi(\Delta)$, the \emph{group of projectivities} of~$\Delta$, which is a combinatorial but not a topological
  invariant of~$\Delta$.  This group is studied for combinatorial manifolds and, in particular, for polytopal simplicial
  spheres.  The results are applied to a coloring problem for simplicial (or, dually, simple) polytopes which arises in
  the area of toric manifolds.
\end{abstract}

\section{Introduction}

In \cite{733.52006} Davis and Januszkiewicz introduce an $(n+d)$-dimensional smooth manifold $\cZ_P$ built from a
$d$-dimensional simple convex polytope~$P$ with $n$~facets.  These manifolds play a significant role in the study of
\hbox{(quasi-)}toric manifolds.  We briefly sketch the construction.  Let $P$ be a simple $d$-polytope with $n$~facets.
Fix an ordering of the facets $\cF=(F_1,\ldots,F_n)$ and let $T$ be the $n$-dimensional complex algebraic torus
$(\CC\setminus\{0\})^\cF$.  On the product $P\times T$ define an equivalence relation~$\sim$, where $(p,s)\sim(q,t)$ if
and only if $p=q$ and the $i$-th component of the quotient $st^{-1}$ in the group~$T$ is trivial for all facets $F_i$
not containing the point~$p=q$.  We obtain a manifold $\cZ_P$ as the quotient space $(P\times T)/\!\sim$.  For a survey
on the subject see Buchstaber and Panov~\cite{math.AT/0010073}, where the construction of the manifold~$\cZ_P$ is
discussed in Section~3.1. The obvious action of the torus~$T$ on~$\cZ_P$ is free over the interior of~$P$.  Points which
are contained in the relative interior of a $k$-dimensional face have a $(d-k)$-dimensional isotropy group.  In
particular, the isotropy group of each vertex has dimension~$d$.  Buchstaber suggested to study quotients of~$\cZ_P$ by
freely acting subgroups of~$T$, see~\cite[Section~4.4]{math.AT/0010073}.  In this context he defines $s(P)$ as the
maximal dimension of a subgroup of~$T$ which acts freely on~$\cZ_P$.  Izmestiev~\cite{ColoredFacets} defines the
chromatic number~$\gamma(P)$ of~$P$ as the minimal number of colors required to color the facets of~$P$ such that any
two facets sharing a vertex have distinct colors.  He shows that $s(P)\ge n-\gamma(P)$,
see~\cite[4.4.5]{math.AT/0010073}, whereas it is clear that $s(P)\le n-d$, see~\cite[4.4.2]{math.AT/0010073}.

From our main result (Theorem~\ref{thm:main}), which is a statement on combinatorial manifolds, we infer a combinatorial
characterization for the simple $d$-polytopes with $\gamma(P)=d$.  The aforementioned results imply that for such
polytopes we have $s(P)=n-d$.  This gives a partial answer to Problem~4.4.1 in~\cite{math.AT/0010073}.  The case
$s(P)=n-d$ seems to be the most interesting one in this context.  The result for simple $3$-polytopes is classical.  The
result for the special case of simple zonotopes is implicit in~\cite[Lemma~4.2.6]{924.53033} of Davis, Januszkiewicz,
and Scott.  Moreover, since the original submission of this paper I learned that Edwards~\cite{Edwards:balanced} had
announced a solution to the coloring problem already in~1977.  However, to the best of my knowledge, no proof was
published.

The paper is organized as follows.  We start by associating a finite group to each facet of a finite-dimensional
simplicial complex, the \emph{group of projectivities}.  For strongly connected complexes the isomorphism class of the
group does not depend on the facet chosen.  In the next section we investigate the groups of projectivities of
combinatorial manifolds.  It turns out that, in order to determine the group of projectivities, it suffices to have
combinatorial information about the fundamental group plus local combinatorial data.  This result is then specialized to
the case of simplicial spheres which arise as boundaries of convex polytopes.  A polytope is simple if and only if the
boundary of its dual is a simplicial sphere.  Hence we can apply our results on combinatorial manifolds to any simple
polytope~$P$.  This way we obtain the desired result: $\gamma(P)=d$ if and only if each $2$-face of~$P$ has an even
number of vertices.  We conclude the paper with a few remarks and a short appendix on how our results are related to
known results in graph theory.

I am indebted to Ivan Izmestiev and Friederike K\"orner for stimulating discussions on the subject.  Thanks to Carsten
Lange, Julian Pfeifle, and G\"unter M.~Ziegler for giving helpful comments on a previous version of this paper.
Moreover, I am grateful to the anonymous referee for bringing to my attention the paper of Davis, Januszkiewicz, and
Scott~\cite{924.53033}.

\section{Simplicial Complexes}

An \emph{(abstract) simplicial complex} on the vertex set~$V$ is a non-empty collection~$\Delta$ of finite subsets
of~$V$, which is closed with respect to forming subsets.  If $\sigma\in\Delta$ with $\#\sigma=k+1$, we say that the
\emph{simplex}~$\sigma$ has \emph{dimension}~$k$, and we write $\dim\sigma=k$.  Define
$\dim\Delta=\sup\SetOf{\dim\sigma}{\sigma\in\Delta}$.  Throughout the rest of the paper we always assume that
$\dim\Delta<\infty$.  A simplex of~$\Delta$ which is maximal with respect to inclusion is called a \emph{facet}.  If a
simplex~$\sigma$ is contained in another simplex~$\tau$, then $\sigma$ is a \emph{face} of~$\tau$.  The complex~$\Delta$
is called \emph{pure} if all its facets have the same dimension.  The maximal proper faces of the facets are the
\emph{ridges}.  For a given face $\sigma\in\Delta$, the \emph{(closed) star} $\st\sigma$ is the subcomplex generated by
the facets containing~$\sigma$, whereas the \emph{link} $\lk\sigma$ is the subcomplex of~$\st\sigma$ of faces not
containing~$\sigma$.

By introducing barycentric coordinates on the simplices and extending according to the concept of weak topology, every
finite-dimensional simplicial complex $\Delta$ defines a locally compact and metrizable Hausdorff space $\norm{\Delta}$,
which is compact if and only if $\Delta$ is finite; see any topology textbook, e.g.\ Munkres~\cite{Munkres84:0}, for the
details.  We frequently apply notions from topology to~$\Delta$ which, if no confusion can arise, are meant to refer
to~$\norm{\Delta}$.

The \emph{dual graph} $\Gamma(\Delta)$ of~$\Delta$ is an abstract graph whose nodes are the facets of~$\Delta$, and
where an edge between two facets corresponds to a common ridge.  We call $\Delta$ \emph{strongly connected} if the graph
$\Gamma(\Delta)$ is connected.  Strong connectedness clearly implies connectedness in the topological sense.  Moreover,
if $\Delta$ is strongly connected, then $\Delta$ is pure.  However, our definition of the dual graph also makes sense
for non-pure complexes.  In the non-pure case each connected component of the dual graph consists of facets of the same
dimension.

For each ridge $\rho$ contained in two facets $\sigma$, $\tau$, there is a unique vertex~$v(\sigma,\tau)$ which is
contained in~$\sigma$ but not in~$\tau$.  We define the \emph{perspectivity} $\langle\sigma,\tau\rangle:\sigma\to\tau$ by setting
$$w\mapsto\left\{\begin{array}{cl}v(\tau,\sigma)&\text{if $w=v(\sigma,\tau)$},\\
    w&\text{otherwise.}\end{array}\right.$$

Let $g=(\sigma_0,\sigma_1,\ldots,\sigma_n)$ be a \emph{facet path} in~$\Gamma(\Delta)$, that is, for each each~$i$ the
facets~$\sigma_i$ and~$\sigma_{i+1}$ share a common ridge.  The \emph{projectivity} $\langle g\rangle$ from $\sigma_1$
to~$\sigma_n$ along~$g$ is the concatenation
$$\langle
g\rangle=\langle\sigma_0,\sigma_1,\ldots,\sigma_n\rangle=\langle\sigma_0,\sigma_1\rangle\langle\sigma_1,\sigma_2\rangle\cdots\langle\sigma_{n-1},\sigma_n\rangle$$
of perspectivities.  The map~$\langle g\rangle$ is a bijection from $\sigma_0$ to~$\sigma_n$.  The facet path~$g$ is
\emph{closed} if $\sigma_0=\sigma_n$.  A closed facet path from~$\sigma_0$ to~$\sigma_0$ is called a \emph{facet loop}
based at~$\sigma_0$.  We denote the concatenation of two facet paths $g=(\sigma_0,\sigma_1,\ldots,\sigma_n)$ and
$h=(\sigma_n,\sigma_{n+1},\ldots,\sigma_m)$ by $g*h$.  Clearly, $\langle g*h\rangle=\langle g\rangle\langle h\rangle$.

For a given facet~$\sigma_0$ the set of projectivities along facet loops based at~$\sigma_0$ forms a
group~$\Pi(\Delta,\sigma_0)$, the \emph{group of projectivities} of~$\Delta$ at~$\sigma_0$.  The group of projectivities
is a (permutation) subgroup of the symmetric group~$\Sym\sigma_0$, the group of all bijections on the set of vertices
of~$\sigma_0$.  The inverse of the facet path~$g$ is denoted by~$g^-$.

\begin{lem}\label{lem:path}
  Let $g$ be a facet path from the facet~$\sigma_0$ to the facet~$\sigma_1$.  Then $\Pi(\Delta,\sigma_0)=\langle
  g\rangle\,\Pi(\Delta,\sigma_1)\,\langle g^-\rangle$.
\end{lem}

This implies that for strongly connected~$\Delta$ the isomorphism class of~$\Pi(\Delta,\sigma)$ does not depend on the
choice of the \emph{base facet}~$\sigma_0$.  We write~$\Pi(\Delta)$ and call it the \emph{group of projectivities}
of~$\Delta$.  In this case the group of projectivities is a combinatorial invariant of~$\Delta$.

Let $\Delta$ and $\Delta'$ be finite-dimensional simplicial complexes and let $f:\Delta\to\Delta'$ a \emph{simplicial
  map}, that is, $f$ is a map between the vertex sets which preserves the inclusion among the faces.  A simplicial map
is called \emph{non-degenerate} if it preserves dimension.  Consider a non-degenerate simplicial map between simplicial
complexes of the same dimension.  In this case facet loops are mapped to facet loops and we obtain an induced map
$$f_\#=\Pi(f):\Pi(\Delta,\sigma_0)\to\Pi(\Delta',f(\sigma_0)).$$

Consider the category $\cC_d$ of pairs $(\Delta,\sigma_0)$, where $\Delta$ is a simplicial complex of fixed
dimension~$d$ and $\sigma_0$ is a facet of~$\Delta$.  As morphisms take the non-degenerate simplicial maps which map
base facets to base facets.
  
\begin{prp}
  $\Pi(\cdot)$ is a covariant functor from the category~$\cC_d$ into the category of finite groups.
\end{prp}

\begin{prp}\label{prp:subcomplex}
  Let $\Delta$ and~$\Delta'$ be $d$-dimensional simplicial complexes, and let $\sigma_0$ be a facet of~$\Delta$.  If
  $f:(\Delta,\sigma_0)\to(\Delta,f(\sigma_0))$ is a non-degenerate simplicial map which is injective if restricted to
  the set of facets, then the induced map~$f_\#$ is a group monomorphism.  In particular, the group of projectivities of 
  a full-dimensional subcomplex of~$\Delta$ which contains~$\sigma_0$ is a subgroup of~$\Pi(\Delta,\sigma_0)$.
\end{prp}

We want to determine the groups of projectivities for simplicial complexes~$\Delta$ of dimension at most~$1$.  Up to an
isomorphism there is a unique simplicial complex of dimension~$-1$, namely $\{\emptyset\}$.  It has a unique facet and
no ridges, so its dual graph consists of a single node.  Its group of projectivities~$\Pi(\{\emptyset\},\emptyset)$ is
trivial.  Similarly, if $\dim\Delta=0$ the facets correspond to the vertices, and $\Pi(\Delta,\{v\})$ is trivial for any
vertex~$v$.  The $1$-dimensional simplicial complexes are precisely the graphs.  The edges are the facets (except for
possibly existing isolated nodes).  Each edge has two nodes, so the group of projectivities is of order at most~$2$.

\begin{prp}\label{prp:n-gon}
  Let $\Delta$ be a graph and let $\sigma_0=\{v,w\}$ be an edge.  If the connected component of~$\sigma_0$ in~$\Delta$
  contains an odd cycle, then $\Pi(\Delta,\sigma_0)$ is generated by the transposition~$(v\ w)$.  Otherwise, the
  connected component of~$\sigma_0$ is bipartite and $\Pi(\Delta,\sigma_0)$ is trivial.
\end{prp}

There are many ways to build new complexes from given ones.  We will explore one construction and its impact on the
group of projectivities.

Let $\Delta$ and $\Delta'$ be finite dimensional simplicial complexes over the vertex sets~$V$ and~$V'$, respectively,
where $V$ is disjoint from~$V'$.  The \emph{join} of $\Delta$ and~$\Delta'$ is defined to be
$$\Delta*\Delta'=\SetOf{\sigma\cup\sigma'}{\sigma\in\Delta,\ \sigma\in\Delta'}.$$
Clearly,
$\dim\Delta*\Delta'=\dim\Delta+\dim\Delta'+1$.  The facets of $\Delta*\Delta'$ are unions of facets of~$\Delta$ and
$\Delta'$; the ridges are unions of a facet of one complex with a ridge of the other.  Forming the join of two complexes
is, in fact, a topological operation: $\norm{\Delta*\Delta'}$ is homeomorphic to the double mapping cylinder of the
projections $\norm{\Delta\times\Delta'}\to\norm{\Delta},\norm{\Delta'}$.

\begin{prp}\label{prp:join}
  Let $\Delta$ and~$\Delta'$ both be finite-dimensional simplicial complexes with facets~$\sigma_0$ and~$\sigma_0'$,
  respectively.  Then $\Pi(\Delta*\Delta',\sigma_0\cup\sigma_0')=\Pi(\Delta,\sigma_0)\times\Pi(\Delta',\sigma_0')$.
\end{prp}

\begin{proof}
  We claim equality instead of the mere existence of an isomorphism because the direct product can be interpreted as an
  inner direct product as follows.  The maps $f:\Delta\to\Delta*\Delta':\sigma\mapsto\sigma\cup\sigma_0'$ and
  $f':\Delta'\to\Delta*\Delta':\sigma'\mapsto\sigma_0\cup\sigma'$ both are non-degenerate and injective, which yields
  monomorphisms $\Pi(f)$ and $\Pi(f')$, respectively, by Proposition~\ref{prp:subcomplex}.
  
  We have to prove that each projectivity in the join can be written as a product of a projectivity in~$\Delta$ with a
  projectivity in~$\Delta'$.
  
  The simplices of~$\Delta*\Delta'$ are written as $\sigma\cup\sigma'$, implying that $\sigma\in\Delta$ and
  $\sigma'\in\Delta'$.  Note that the distinct facets~$\sigma\cup\sigma'$ and~$\tau\cup\tau'$ are adjacent if and only
  if $\sigma=\tau$ and $\sigma'$ adjacent to~$\tau'$ in~$\Delta'$, or $\sigma$ adjacent to~$\tau$ in~$\Delta$ and
  $\sigma'=\tau'$.
  
  Moreover, any two facets~$\sigma\cup\sigma'$ and~$\tau\cup\tau'$ with $\sigma$ adjacent to~$\tau$ and $\sigma'$
  adjacent to~$\tau'$ are contained in the star of the codimension-$2$-face $(\sigma\cap\tau)\cup(\sigma'\cap\tau')$.
  There are precisely two more facets contained in this star, namely $\sigma\cup\tau'$ and~$\tau\cup\sigma'$.  Applying
  the Propositions~\ref{prp:subcomplex} and~\ref{prp:n-gon} to $\st((\sigma\cap\tau)\cup(\sigma'\cap\tau'))$ yields
  $$\langle\sigma\cup\sigma',\sigma\cup\tau',\tau\cup\tau',\tau\cup\sigma',\sigma\cup\sigma'\rangle = 1$$
  and thus
  \begin{equation}\label{eq:product}
  \langle\sigma\cup\sigma',\sigma\cup\tau',\tau\cup\tau'\rangle = \langle\sigma\cup\sigma',\tau\cup\sigma',\tau\cup\tau'\rangle.
  \end{equation}
  
  Invoking the identity~(\ref{eq:product}) several times, allows to ``sort'' a projectivity: Each
  projectivity~$\pi$ from~$\sigma_0\cup\sigma_0'$ onto itself can be written as the product
  \begin{eqnarray*}
    \pi
    &=&\langle\sigma_0\cup\sigma_0',\sigma_1\cup\sigma_0',\ldots,(\sigma_m=\sigma_0)\cup\sigma_0')\rangle\\
    &&\langle\sigma_0\cup\sigma_0',\sigma_0\cup\sigma_1',\ldots,\sigma_0\cup(\sigma_n'=\sigma_0')\rangle.
  \end{eqnarray*}
\end{proof}\goodbreak

Izmestiev~\cite{IvanPrivateComm} has proved a partial converse of the previous proposition.

The $d$-dimensional simplicial complex~$\Delta$ on the vertex set~$V$ is called \emph{balanced} if there is a map
$c:V\to\{0,\ldots,d\}$ such that whenever $\{v,w\}$ is an edge in~$\Delta$ then $c(v)\ne c(w)$.  The map $c$ is called a
\emph{proper $d$-coloring} of~$\Delta$.  Clearly, a proper $d$-coloring of~$\Delta$ is the same as a simplicial
projection from~$\Delta$ onto the standard $d$-simplex which is injective on each simplex.  Occasionally, this is called
a \emph{folding map} of~$\Delta$.  Important examples for balanced simplicial complexes are provided by Coxeter
complexes and Tits buildings, see Stanley~\cite[pp~104ff]{838.13008}.  For properties of proper colorings or folding
maps in the context of toric manifolds see Davis and Januszkiewicz~\cite[Lemma~1.14 and Example~1.15]{733.52006}.

We call a simplicial complex \emph{locally strongly connected} if it is strongly connected and, additionally, the star
of each vertex is also strongly connected.  There are strongly connected complexes which is not locally strongly
connected.  For instance, consider a $2$-dimensional complex whose dual graph is a path such that the two triangles
corresponding to the end points of the path share a unique vertex~$v$.  The star of~$v$ is not strongly connected.

\begin{prp}\label{prp:balanced}
  Let $\Delta$ be a locally strongly connected simplicial complex.  Then $\Delta$ is balanced if and only if
  $\Pi(\Delta)$ is trivial.
\end{prp}

\begin{proof}
  Fix an arbitrary facet $\sigma_0$ of~$\Delta$ and an arbitrary coloring of the vertices of~$\sigma_0$.  For each facet
  path~$g$ from $\sigma_0$ to some other facet~$\sigma$ the projectivity~$\langle g\rangle$ induces a coloring of the
  vertices of~$\sigma$.  Two such colorings induced by facet paths $g$ and~$g'$, respectively, coincide if and only if
  the projectivity $\langle g'*g^-\rangle=\langle g'\rangle\langle g\rangle^{-1}$, which is induced by the facet loop
  $g'*g^-$ based at $\sigma_0$, is the identity.  Observe that, in general, the color of a vertex~$v\in\sigma$ does
  depend on the choice of the facet~$\sigma$.  Since, however, the star of~$v$ in~$\Delta$ is also strongly connected,
  this color is the same for all facets containing~$v$.\\ \nothing\hfill
\end{proof}

It is worth mentioning that the property of being balanced is by no means a topological invariant.  To the contrary, for 
arbitrary~$\Delta$ the barycentric subdivision $\sd\Delta$ is always balanced.

\section{Combinatorial Manifolds}

We now impose severe topological restrictions on the simplicial complexes studied.  A finite $d$-dimensional simplicial
complex~$\Delta$ is a \emph{combinatorial manifold} if the link of each $k$-face is a simplicial sphere of dimension
$d-k-1$.  In particular, the link of each codimension-$2$-face is a $1$-sphere, that is, the boundary of a polygon on
the combinatorial level.  Note that our (combinatorial) manifolds are always compact and without boundary.  However, the
results below can suitably be extended to combinatorial manifolds with boundary.

If $\Delta$ is a combinatorial manifold, then $\norm{\Delta}$ is a PL-manifold.  Conversely, a PL-manifold~$M$ always
admits a triangulation~$\Delta$ (compatible with the PL-structure), such that $\Delta$ is a combinatorial manifold.  For
a general introduction to combinatorial and PL-manifolds see Hudson~\cite{189.54507}, Glaser~\cite{212.55603},
or~\cite[65 (IX.17)]{Ito1987}.

Throughout the following let $\Delta$ be a combinatorial manifold.  This implies that the dual graph~$\Gamma(\Delta)$ is
strongly connected, so the isomorphism class of the group of projectivities does not depend on the facet chosen.

Consider the joint geometric realization~$\norm{\Delta}$ of $\Delta$ and its dual block complex~$\Delta^*$ within a
realization of the first barycentric subdivision~$\sd\Delta$, see Munkres~\cite[\S64]{Munkres84:0} and also
Glaser~\cite[pp.~83ff]{212.55603}.  This way each facet path canonically yields an edge path in the $1$-skeleton of the
dual block complex~$\Delta^*$ and vice versa.  Often we will not distinguish between a facet path and its corresponding
edge path in~$\Delta^*$.  As $\Delta$ is a combinatorial manifold the blocks in~$\Delta^*$ are, in fact, cells.
In particular, the blocks are simply connected.

It is known that any path in~$\norm{\Delta}=\norm{\Delta^*}$ is homotopic to a path in the $1$-skeleton of~$\Delta^*$
which is the same as the dual graph of~$\Delta$.  In Seifert and Threlfall~\cite[\S44]{469.55001} this is proved for
simplicial complexes, but the arguments given can directly be extended to arbitrary cell complexes.  In particular, the
fundamental group $\pi_1(\Delta,x_0)$ for $x_0\in\norm\Delta$ is generated by facet loops based at~$\sigma_0$ where
$\sigma_0$ is some facet with $x_0\in\norm{\sigma_0}$.  Usually, in the geometric realization we choose $x_0$ to be the
barycenter of the facet~$\sigma_0$, and we write $\pi_1(\Delta,\sigma_0)$.  Note that, as $\Delta$ is assumed to be
finite, the group $\pi_1(\Delta,\sigma_0)$ is finitely generated.

\def\rPi{{\Pi_0}}

Define the \emph{reduced group of projectivities}~$\rPi(\Delta,\sigma_0)$ to be the subgroup of~$\Pi(\Delta,\sigma_0)$
generated by facet loops based at~$\sigma_0$ which are null-homotopic.  Similar to what is expressed in
Lemma~\ref{lem:path} the reduced group of projectivities is a combinatorial invariant of the connected component
of~$\sigma_0$ in~$\Delta$.

\begin{prp}\label{prp:gen}
  Let $p_1,\ldots,p_m$ be a set of facet loops based at $\sigma_0$ generating the fundamental
  group~$\pi(\Delta,\sigma_0)$.  Then $\Pi(\Delta,\sigma_0)$ is generated by $\rPi(\Delta,\sigma_0)$ together with
  $\langle p_1\rangle,\ldots,\langle p_m\rangle$.
\end{prp}

In particular, if $\pi_1(\Delta,\sigma_0)$ is trivial then $\rPi(\Delta,\sigma_0)=\Pi(\Delta,\sigma_0)$.  The converse
does not hold.

\begin{figure}[htbp]
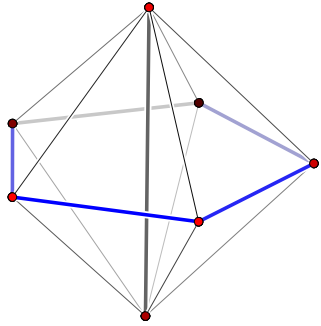

  \begin{center}
    \altgraphicsWidth{file=Lk-kappa.eps}{Lk-kappa.pdf}{6cm}
    \caption{Pentagonal link of an edge in a $3$-manifold.}
    \label{fig:Lk-kappa}
  \end{center}
\end{figure}

The link of each codimension-$2$-face~$\kappa$ is an $n$-gon for some $n\ge 3$; see Figure~\ref{fig:Lk-kappa}.  Due to
the obvious bijection between the facets in~$\lk\kappa$ and the facets in~$\st\kappa$ we see that $\Gamma(\st\kappa)$ is
also an $n$-gon.  The \emph{parity} of~$\kappa$, that is, the property of being \emph{even} or \emph{odd}, is the parity
of~$n$.

Let $\kappa$ be a codimension-$2$-face~$\kappa$, $\sigma$ a facet containing~$\kappa$, and $g$ a path from $\sigma_0$
to~$\sigma$.  As $\st\kappa$ is simply connected we infer that the path $g*l*g^-$ is null-homotopic for any facet
loop~$l$ in~$\st\kappa$ based at~$\sigma$.  Thus we have $\langle g*l*g^-\rangle\in\rPi(\Delta,\sigma_0)$.

If $\kappa$ is odd, in view of Proposition~\ref{prp:n-gon}, the group $\Pi(\st\kappa,\sigma)$ is of order~$2$, generated
by some facet loop~$l$ based at~$\sigma$.  Then $\langle g*l*g^-\rangle$ is a transposition on the set~$\sigma_0$.

\begin{thm}\label{thm:main}
  The reduced group of projectivities $\rPi(\Delta,\sigma_0)$ is generated by the set of all projectivities $\langle
  g*l*g^-\rangle$ where $g$ is a facet path from $\sigma_0$ to some facet~$\sigma$ which contains an odd
  codimension-$2$-face $\kappa$ and $l$ is a facet loop based at~$\sigma$ generating $\Pi(\st\kappa,\sigma)$.  In
  particular, $\rPi(\Delta,\sigma_0)$ is generated by transpositions.
\end{thm}

\begin{proof}
  Let $r$ be an arbitrary facet loop based at~$\sigma_0$ which is null-homotopic.  Without loss of generality let $x_0$
  be the vertex of~$\Delta^*$ corresponding to the barycenter of~$\sigma_0$.  It is known that $r$ can be contracted to
  the constant map $c_{x_0}$ at~$x_0$ within the $2$-skeleton of~$\Delta^*$.  Discretizing a suitable homotopy from $r$
  to $c_{x_0}$ yields a sequence $r_1,\ldots,r_n$ of closed paths in the $1$-skeleton from~$x_0$ to~$x_0$ in the
  $1$-skeleton of~$\Delta^*$ such that $r_1=r$, $r_n=c_{x_0}$, and $r_i$ coincides with $r_{i+1}$ outside some
  $2$-face~$F_i$ of~$\Delta^*$; see Figure~\ref{fig:contract}.  The dual of~$F_i$ in~$\Delta$ is a codimension-$2$-face
  $\kappa_i$.

  \begin{figure}[htbp]
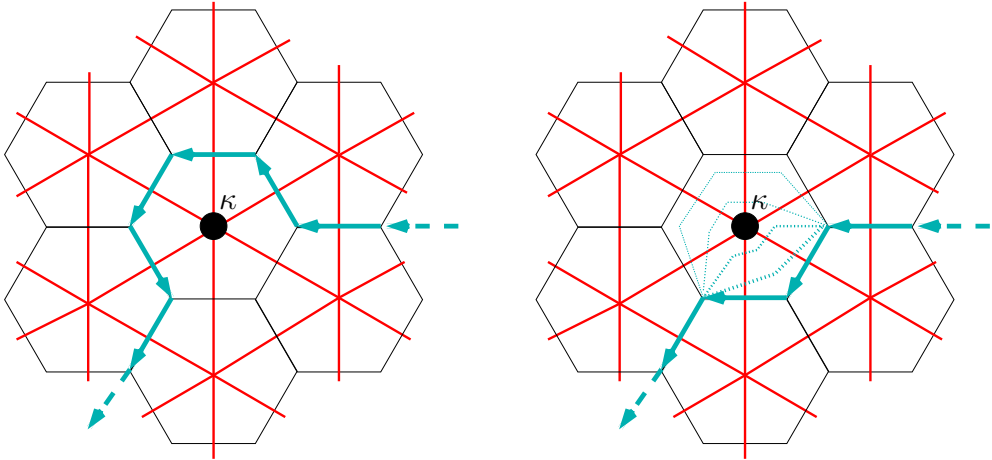

    \begin{center}
      \ifx\pdfoutput\undefined
      \input contract-a.pstex_t\qquad \input contract-b.pstex_t
      \else \includegraphics[angle=270]{contract-a.pdf}\qquad \includegraphics[angle=270]{contract-b.pdf}
      \fi
      \caption{Combinatorial homotopy between paths in the dual graph within the star of
        the codimension-$2$-face~$\kappa$.}
      \label{fig:contract}
    \end{center}
  \end{figure}
  
  Because the facet paths~$r_1$ and~$r_2$ are the same outside $\st\kappa$ we have that the projectivity $\langle
  r_1\rangle\langle r_2\rangle^{-1}$ coincides with some projectivity $\langle g*l*g^-\rangle$, where $g$ is the common
  initial segment of~$r_1$ and $r_2$ up to some facet~$\sigma\in\st\kappa$ and $l$ is a facet loop in~$\st\kappa$ based
  at~$\sigma$.  In particular, by Proposition~\ref{prp:n-gon}, $\langle r_1\rangle\langle r_2\rangle^{-1}$ is either a
  transposition or trivial, depending on the parity of~$\kappa$.
  
  An induction on~$n$ establishes the theorem.
\end{proof}

\begin{cor}\label{cor:class}
  The reduced group of projectivities $\rPi$ of a combinatorial manifold is isomorphic to a direct product of symmetric
  groups.
\end{cor}

\begin{figure}[htbp]
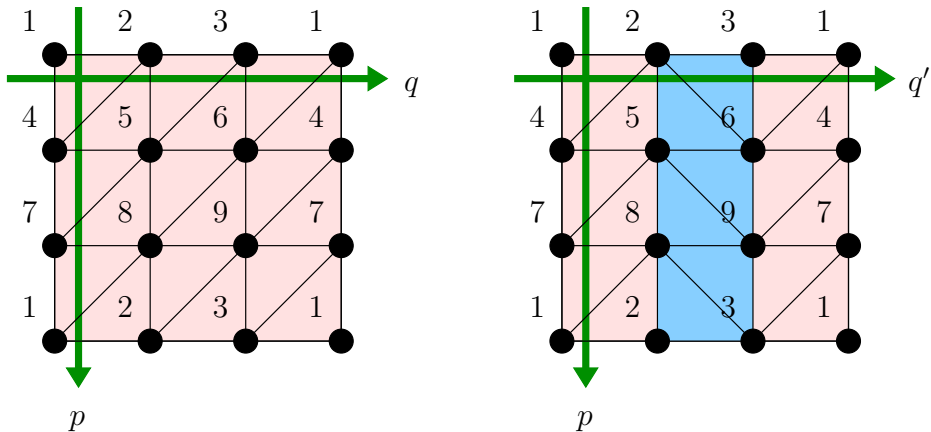

  \begin{center}
    \ifx\pdfoutput\undefined
    \input torus.pstex_t\qquad \input anti-torus.pstex_t
    \else \includegraphics[width=5cm]{torus.pdf}\qquad \includegraphics[width=5cm]{anti-torus.pdf}
    \fi

    \caption{Left: Standard torus~$T$ with facet loops~$p$ and~$q$ which correspond to a generating
      system of the fundamental group. Right: ``anti-torus''~$A$ with similar facet loops~$p$ and~$q'$.}
    \label{fig:torus}
  \end{center}
\end{figure}

The same result does not hold for the whole group of projectivities~$\Pi$.  For an example see Figure~\ref{fig:torus}.

\begin{cor}\label{cor:even}
  The reduced group of projectivities $\rPi(\Delta,\sigma_0)$ is trivial if and only if each codimension-$2$-face
  of~$\Delta$ is even.
\end{cor}

\begin{cor}\label{cor:balanced}
  Suppose that $\Delta$ is simply connected.  Then $\Delta$ is balanced if and only if each codimension-$2$-face
  of~$\Delta$ is even.
\end{cor}

Corollary~\ref{cor:balanced} seems to be known: It is announced, without a proof, in Edwards~\cite{Edwards:balanced}:
``The above theorem [on a reformulation of the Four Color Problem] developed from a lunch table conversation at
I.H.E.S., Bures-sur-Yvette, France, in which P. Deligne-R. MacPherson-J. Morgan observed that a closed, $1$-connected,
PL triangulated $n$-manifold is $(n+1)$-colorable $\Longleftrightarrow$ each $(n-2)$-simplex has even order.''

The group of projectivities is an interesting invariant of a combinatorial manifold.  Consider, for example, two
different triangulations of the $2$-torus~$\Sph^1\times\Sph^1$ as depicted in Figure~\ref{fig:torus}.  The first
triangulation~$T$ (to the left) is standard.  The second triangulation~$A$ is produced from~$T$ by flipping the diagonal
edges in the three squares of the middle column; in order to give it some name, call it \emph{anti-torus}.  Several
combinatorial invariants of~$T$ and~$A$ coincide: e.g., the $f$-vector, the vector of vertex-degrees in the graph, the
Altshuler determinant.  But the groups of projectivities differ.

This can be seen as follows.  Fix the facet $\sigma_0=\{1,2,4\}$ in both triangulations. The codimension-$2$-faces are
the vertices.  In both triangulations, the link of each vertex is the boundary of a hexagon.  Therefore, by
Corollary~\ref{cor:even}, the only potentially non-trivial contributions to $\Pi(T,\sigma_0)$ and $\Pi(A,\sigma_0)$ can
come from the fundamental group which is known to be isomorphic to $\ZZ\times\ZZ$.  The facet loops
$$p=(\{1,2,4\},\{2,4,5\},\{4,5,7\},\{5,7,8\},\{1,7,8\},\{1,2,8\},\{1,2,4\})$$
and
$$q=(\{1,2,4\},\{2,4,5\},\{2,3,5\},\{3,5,6\},\{1,3,6\},\{1,4,6\},\{1,2,4\})$$
generate the group~$\pi_1(T,\sigma_0)$.
Verify that both $\langle p\rangle$ and $\langle q\rangle$ are equal to the identity.  Now $p$ and
$$q'=(\{1,2,4\},\{2,4,5\},\{2,5,6\},\{2,3,6\},\{1,3,6\},\{1,4,6\},\{1,2,4\})$$
generate $\pi_1(A,\sigma_0)$.  Again
$\langle p\rangle=\id$, but $\langle q'\rangle$ is the $3$-cycle $(1\ 4\ 2)$.  Therefore, $\Pi(A,\sigma_0)\cong\ZZ/3$.

\section{Polytopes}

A polytope is \emph{simple} if each of its vertex figures is a simplex, or, equivalently, for any given vertex~$v$ there
is a 1--1 correspondence between the sets of edges through~$v$ and the faces containing~$v$.  For an introduction to
the theory of convex polytopes, see Ziegler~\cite{GMZ}.  Here we restrict our attention to polytopes which are convex.

There is another way to characterize simple polytopes, which suits our needs: A polytope~$P$ is simple if and only if
its dual~$P^*$ is simplicial, that is, each proper face is a simplex.  In particular, the boundary complex of a simple
polytope is the dual cell complex of a polytopal sphere.  Therefore, we can dualize our definition of perspectivity.
The results of the previous section apply.

Let $v$ be a vertex in the simple $d$-polytope~$P$.  Denote the set of facets through~$v$ by $\cF(v)$.  If $w$ is a
vertex adjacent to~$v$ then there is a unique facet $F(v,w)$ contained in~$\cF(v)\setminus\cF(w)$.  The
\emph{perspectivity} from $v$ to~$w$ is defined as
$$\langle v,w\rangle:\cF(v)\to\cF(w):F\mapsto\left\{\begin{array}{cl}F(w,v)&\text{if $F=F(v,w)$,}\\
    F&\text{otherwise.}\end{array}\right.$$
Again \emph{projectivities} are concatenations of perspectivities.  As the
boundary of a polytope is connected the isomorphism class of the group of projectivities does not depend on the vertex
chosen.

Note that each $2$-face of~$P$ corresponds to (the link of) a codimension-$2$-face of the dual.  Therefore the following
corollary follows from our Theorem~\ref{thm:main}.

\begin{cor}\label{cor:simple-vanish}
  For any vertex~$v$ the group of projectivities~$\Pi(P,v)$ is generated by projectivities with respect to paths around
  the $2$-faces with an odd number of vertices.  In particular, if each $2$-face has an even number of vertices, then
  the group of projectivities vanishes.
\end{cor}

\begin{proof}
  The boundary complex of a polytope is homeomorphic to a sphere, and thus the fundamental group is trivial, provided
  that the dimension of the polytope is at least~$3$.  The group of projectivities coincides with the reduced group of
  projectivities.  For $2$-dimensional polytopes the only $2$-face is the polytope itself and the result follows from
  Proposition~\ref{prp:n-gon}.  A $1$-dimensional polytope does not have any $2$-face, its dual graph consists of two
  isolated points, and hence the group of projectivities is trivial.
\end{proof}

This directly allows to compute the group of projectivities of many known polytopes, including all regular simple
polytopes.

\begin{cor}\label{cor:odd-regular}
  The group of projectivities of the $d$-simplex is isomorphic to~$S_d$.
  
  The group of projectivities of the dodecahedron is isomorphic to~$S_3$.
  
  The group of projectivities of the regular $120$-cell is isomorphic to~$S_4$.

  The group of projectivities of the $d$-cube is trivial.
\end{cor}

\begin{proof}
  Each $2$-face of a simplex is a triangle.  Each $2$-face of the dodecahedron and the $120$-cell is a pentagon.  Each
  $2$-face of the $d$-cube is a quadrangle.
\end{proof}

In Proposition~\ref{prp:join} we discussed the effect of forming joins of simplicial complexes on the group of
projectivities.  This can be translated into a result about simple polytopes.

\begin{cor}\label{cor:product}
  Let $P$ and $Q$ be simple polytopes with respective vertices $v$ and~$w$.  Then $\Pi(P\times
  Q,(v,w))=\Pi(P,v)\times\Pi(Q,w)$.
\end{cor}

\begin{proof}
  The product $P\times Q$ is again a simple polytope.  Its boundary complex is dual to the join of the duals of the
  boundary complexes of~$P$ and~$Q$.
\end{proof}

The example of products of simplices shows that for any \emph{partition} of~$d$, that is, a sequence $(d_1,\ldots,d_k)$
of natural numbers with $d_i\ge 1$ and $\sum d_i=d$, there is a simple $d$-polytope whose group of projectivities is
isomorphic to~$S_{d_1}\times\cdots\times S_{d_k}$.  From Corollary~\ref{cor:class} we infer that, in fact, this is the
only class of groups which occurs as groups of projectivities of simple polytopes.  We obtain a combinatorial invariant
of a simple polytope.

\begin{cor}
  Let $P$ be a simple $d$-polytope.  Then there is a unique partition~$(d_1,\ldots,d_k)$ of~$d$ with $d_1\le
  d_2\le\cdots\le d_k$ such that $\Pi(P)\cong S_{d_1}\times\cdots\times S_{d_k}$.
\end{cor}

\begin{figure}[htp]
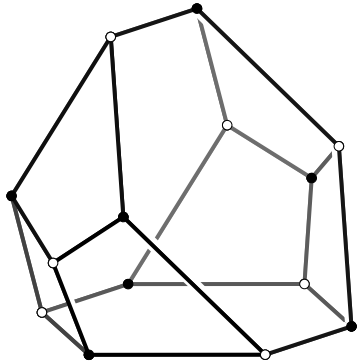

  \begin{center}
    \altgraphicsWidth{file=minimal-even3-corrected.eps}{minimal-even3-corrected.pdf}{6cm}
    \caption{Even simple $3$-polytope~$M$ which is not combinatorially equivalent to any zonotope.  The polytope~$M$ is
      constructed as the \emph{blending} of two cubes. This operation on simple polytopes has been introduced by
      Barnette~\cite{174.25404} as \emph{joining}; it is called \emph{connected sum} in~\cite{math.AT/0010025}.  The
      $f$-vector of~$M$ equals $(14,21,9)$.  This is not the $f$-vector of any zonotope, because zonotopes, being
      centrally symmetric, have an even number of facets.  One can show that among all even simple $3$-polytopes which
      are not combinatorially equivalent to any zonotope, $P$ has the minimal number of vertices as well as the minimal
      number of facets.  We indicate the bipartition of the vertex set according to Theorem~\ref{thm:coloring}.  The
      picture has been produced with polymake~\cite{PolymakeCite} and JavaView~\cite{javaview}.}
    \label{fig:even}
  \end{center}
\end{figure}

The Corollary~\ref{cor:simple-vanish} characterizes those simple polytopes whose $2$-faces have an even number of
vertices.  We call such simple polytopes \emph{even}.  Note that each simple zonotope is an even simple polytope.  But,
an easy construction shows that the even simple polytopes form a (much) wider class.  For an example see
Figure~\ref{fig:even}.

\def\overGamma{{\overline\Gamma}}

Let $P$ be an arbitrary $d$-polytope. Define a graph~$\overGamma(P)$ whose nodes are the facets of~$P$; two facets are
joined by an edge in~$\overGamma(P)$ if their intersection is not empty.  A \emph{proper (node) coloring} of a graph is
an assignment of a color to each node such that any two adjacent nodes have different colors.  The \emph{chromatic
  number} of a graph is the minimal number of colors in a proper coloring.  Following Izmestiev~\cite{ColoredFacets},
the \emph{chromatic number}~$\gamma(P)$ of the polytope~$P$ is now defined as the chromatic number of the
graph~$\overGamma(P)$.  As every vertex of~$P$ is contained in at least~$d$ facets, it is clear that it requires at
least~$d$ colors to color~$\overGamma(P)$ properly.  Moreover, if $\gamma(P)=d$ then $P$ is simple.

The $1$-skeleton of a polytope also forms an abstract graph, which is more commonly studied in polytope theory.  In
order to avoid confusion we call this graph the \emph{vertex-edge-graph} of~$P$.

For simple polytopes the graph~$\overGamma(P)$ coincides with its dual graph, that is, the vertex-edge-graph of the dual
(simplicial) polytope: This follows from the fact that each vertex figure of a simple polytope is a simplex.  Hence any
two facets which share a vertex already have a common ridge.
\goodbreak

\begin{thm}\label{thm:coloring}
  Let $P$ be a simple $d$-polytope.  Then the following properties are equivalent.
  \begin{enumerate}
  \item The polytope~$P$ is even.
  \item The vertex-edge-graph of~$P$ is bipartite.
  \item The boundary complex~$\partial P^*$ of the dual is balanced.
  \item $\gamma(P)=d$.
  \end{enumerate}
\end{thm}

\begin{proof}
  Let $P$ be an even simple $d$-polytope.  Due to Corollary~\ref{cor:simple-vanish} we know that this is characterized
  by the property that the group of projectivities vanishes.  A proper coloring of the facets of~$P$ clearly corresponds
  to a proper coloring of the vertices of the dual~$P^*$.  The existence of such a coloring now follows from
  Proposition~\ref{prp:balanced}.  This proves the equivalence of the first, the third and the fourth statement.  The
  equivalence of the first and the second statement is known.  We indicate a short proof in the Appendix.
\end{proof}

The same result for $3$-dimensional polytopes is classical, see Ore~\cite[13.1.1]{149.21101} and also
Izmestiev~\cite{ColoredFacets} for a more recent proof.  The proofs employ techniques, for which it seems to be unclear
how they can be generalized to higher dimensions.  The result for $4$-dimensional polytopes follows from work of Goodman
and Onishi~\cite{397.05021}.  Davis, Januszkiewicz, and Scott proved in~\cite[Lemma~4.2.6]{924.53033} that the boundary
complex of the dual of a simple zonotope is balanced.

Recall the definition of the manifold $\cZ_P=(P\times T)/\!\sim$ from the introduction.  The number $s(P)$ is defined as
the maximal dimension of a subgroup of the algebraic torus~$T$ which acts freely on~$\cZ_P$.

\begin{cor}
  If $P$ is an even simple $d$-polytope, then $s(P)=n-d$.
\end{cor}

\begin{proof}
  The dimension $s(P)$ of a freely acting subgroup is bounded from above by $n-d$ according to Buchstaber and
  Panov~\cite[4.4.2]{math.AT/0010073}.  The same number is bounded from below by $n-\gamma(P)$ by a result of
  Izmestiev~\cite{ColoredFacets}, see~\cite[4.4.5]{math.AT/0010073}.  But the Theorem~\ref{thm:coloring} enforces
  $\gamma(P)=d$.
\end{proof}

From the fact that the coloring defined in the proof of Proposition~\ref{prp:balanced} is indeed a proper coloring we
immediately obtain the following corollary.

\begin{cor}
  Let $P$ be an even simple $d$-polytope. Suppose that $v$ and $w$ are adjacent vertices of~$P$.  Then the two facets
  $F(v,w)$ and $F(w,v)$ are disjoint.
\end{cor}

\section{Concluding Remarks}

The terms \emph{perspectivity} and \emph{projectivity} are borrowed from incidence geometry, in particular from the
theory of projective planes and generalized polygons, see Van Maldeghem~\cite[Section~1.5]{PolygonBook}.  These notions
in turn are inspired by concepts from projective geometry.  Moreover, some properties of our groups of projectivities
suggest that they can also be seen as some combinatorial analogue of holonomy groups.

It is natural to ask what kind of finite groups can arise as the groups of projectivities of interesting simplicial
complexes.  From Theorem~\ref{thm:main} we know that the group of projectivities of any simply connected combinatorial
manifold is necessarily isomorphic to a, possibly trivial, product of symmetric groups.
Izmestiev~\cite{IvanPrivateComm} shows that for each conjugacy class of a subgroup of the symmetric group~$S_{d+1}$ of
degree~$d+1$ there is a combinatorial manifold such that the given group arises as the group of projectivities.

A lot is known about the $f$-vectors of balanced simplicial complexes.  This is particularly true for balanced
Cohen-Macaulay complexes which include the boundary complexes of simplicial polytopes.  See
Stanley~\cite[Section~III.4]{838.13008} as well as Billera and Bj\"orner~\cite[15.1.3, 15.2.4]{DCG:FaceNumbers}.

There is an intriguing question on planar graphs which is open for quite some time now.  It might be worthwhile to
explore whether the methods developed in this paper can contribute towards a solution.

\begin{con}{\rm (Barnette 1970)}
  The vertex-edge-graph of an even simple $3$-polytope contains a Hamiltonian cycle.
\end{con}

For one fairly large class of even simple polytopes one can see immediately that this conjecture holds.  Start from an
arbitrary simple polytope~$P$.  Successively truncate all the faces with increasing dimension to obtain~$P'$.
Truncation is dual to stellar subdivision.  So the boundary of the polytope~$P'$ is dual to the barycentric subdivision
of the boundary of $P^*$.  In particular, $P'$ is even.  Now, each spanning tree in the dual graph of~$P$ yields a
Hamiltonian cycle in the vertex-edge-graph of~$P'$.

\section{Appendix}

Let $\Gamma$ be a finite graph with node set~$V$ and edge set~$E$.  Consider the $\GF{2}$-vector space $\GF{2}^E$ of
mappings of $E$ into~$\GF{2}$. Each subset of~$E$ corresponds to such a map via the characteristic function.  The
\emph{cycle space} of~$\Gamma$ is the subspace~$C(\Gamma)$ of~$\GF{2}^E$ generated by all cycles of~$\Gamma$.

A \emph{pure polytopal complex} is a finite collection $P_1,\ldots,P_k\subset\RR^n$ of convex $d$-polytopes such that
the intersection of any two polytopes is a face in both.  The boundary complex of any polytope is a polytopal complex,
for instance.  We want to recursively define the \emph{constructibility} of a polytopal complex: A polytope is
constructible.  A pure polytopal complex~$\Delta$ which is the union of pure constructible subcomplexes $A$ and~$B$ is
constructible if the intersection $A\cap B$ is a pure constructible complex.  The notion of constructibility generalizes
the concept of \emph{shellability}, see Ziegler~\cite[\S8]{GMZ}.  From a theorem of Bruggesser and Mani~\cite{251.52013}
it is known that the boundary complexes of polytopes are shellable and thus constructible.

The $1$-skeleton of a polytopal complex forms an abstract graph~$\Gamma(\Delta)$.  For $\Delta$ being the boundary of a
convex polytope we called $\Gamma(\Delta)$ the \emph{vertex-edge-graph} of the polytope above.  The following result is
known.  A proof follows from a double induction on the dimension of the complex~$\Delta$ and the number of the polytopes
comprising~$\Delta$.

\begin{prp}
  Let $\Delta$ be a constructible polytopal complex.  Then the cycle space $C(\Gamma(\Delta))$ is generated by the
  cycles corresponding to the $2$-faces of~$\Delta$.
\end{prp}

A finite graph is bipartite if and only if all the cycles in a cycle basis have even length.  In particular, a simple
polytope is even if and only if its graph is bipartite.  This proves the equivalence of the first and the third
statement in Theorem~\ref{thm:coloring}.

The vertex-edge-graph of any simple $d$-polytope is $d$-regular.  A bipartite regular graph has an even number of
vertices because, by double counting, both color classes are of the same size.

\begin{cor}
  An even simple polytope has an even number of vertices.
\end{cor}

We want to explore the relationship between proper facet colorings of a simple polytope and proper edge colorings of its
vertex-edge-graph.  An edge coloring of a graph is \emph{proper} if any two edges which share a vertex have distinct
colors.

\begin{prp}\label{prp:edge-coloring}
  Let $P$ be a simple $d$-polytope, and let $c$ be a proper coloring of~$\overGamma(P)$ with $d$~colors.  Then $c$
  induces a proper edge coloring of the vertex-edge-graph~$\Gamma(P)$ with $d$~colors.
\end{prp}

\begin{proof}
  Let $e=\{v,w\}$ be an edge of~$P$.  If $\overGamma(P)$ is properly $d$-colored, then the two facets $F(v,w)$
  and~$F(w,v)$ have the same color.  Assign this color to the edge~$e$.  Evidently, this procedure requires exactly
  $d$~colors.  Assume that this edge coloring is not proper, that is, there are vertices $u$, $v$, $w$ such that
  $\{u,v\}$ and $\{v,w\}$ are edges of the same color.  Then we have $c(F(v,u))=c(F(v,w))$, but the facets $F(v,u)$
  and~$F(v,w)$ both contain the vertex~$v$.  This contradicts the assumption that $c$ is a proper coloring of the
  facets.
\end{proof}

As already mentioned, the graph~$\Gamma(P)$ of~$P$ is $d$-regular.  By a result of Vizing and Gupta, see
West~\cite[6.1.7]{845.05001}, the edges of~$\Gamma(P)$ can be properly colored with at most $d+1$~colors.  K\"onig
proved that a bipartite $d$-regular graph is edge $d$-colorable, see West~\cite[6.1.5]{845.05001}.  Therefore, in view
of Theorem~\ref{thm:coloring}, Proposition~\ref{prp:edge-coloring} can be interpreted as a very special instance of a
classical result from graph theory.

\providecommand{\bysame}{\leavevmode\hbox to3em{\hrulefill}\thinspace}

\goodbreak

\noindent
Michael Joswig\\
Technische Universit\"at Berlin\\
Fakult\"at II: Mathematik und Naturwissenschaften\\
Institut f\"ur Mathematik, MA 6-2\\
Stra\ss{}e des 17.~Juni 136\\
D-10623 Berlin, Germany\\
\texttt{joswig@math.tu-berlin.de}
\end{document}